\documentclass[titlepage,draft,12pt]{article} 
\usepackage[title]{appendix}
\usepackage{amssymb,amsthm} 
\usepackage[a4paper]{geometry}
\usepackage[usenames,dvipsnames]{xcolor}

\date{}

\newcommand{\np}{N^{\perp}}

\newcommand{\ep}{\varepsilon}
\renewcommand{\qed}{{\penalty 10000\mbox{$\quad\Box$}}}
\newcommand{\re}{\mathbb{R}}

\newcommand{\n}{\mathbb{N}}

\newtheorem{thm}{Theorem}[section]
\newtheorem{thmbibl}{Theorem}

\newtheorem{rmk}[thm]{Remark}

\newtheorem{lemma}[thm]{Lemma}

\title{{A concrete realization of the slow-fast alternative for a semi
linear heat equation with homogeneous Neumann boundary conditions}}

\author{Marina Ghisi\vspace{1ex}\\ 
{\normalsize Universit\`a degli Studi di Pisa} \\
{\normalsize Dipartimento di Matematica}\\ 
{\normalsize PISA (Italy)}\\
{\normalsize e-mail: \texttt{marina.ghisi@unipi.it}}
\and
Massimo Gobbino\vspace{1ex}\\ 
{\normalsize Universit\`a degli Studi di Pisa} \\
{\normalsize Dipartimento di Matematica}\\ 
{\normalsize PISA (Italy)}\\  
{\normalsize e-mail: \texttt{massimo.gobbino@unipi.it}}
\and
Alain Haraux\vspace{1ex}\\ 
{\normalsize Universit\'{e} Pierre et Marie Curie} \\
{\normalsize Laboratoire Jacques-Louis Lions}\\ 
{\normalsize PARIS (France)}\\  
{\normalsize e-mail: \texttt{haraux@ann.jussieu.fr}}}

\begin{document}
\maketitle
\begin{abstract}
	We investigate the asymptotic behavior of solutions to a
	semilinear heat equation with homogeneous Neumann boundary
	conditions.
	
	It was recently shown that the nontrivial kernel of the linear
	part leads to the coexistence of fast solutions decaying to 0
	exponentially (as time goes to infinity), and slow solutions
	decaying to 0 as negative powers of $t$.
	
	Here we provide a characterization of slow/fast solutions in terms
	of their sign, and we show that the set of initial data giving
	rise to fast solutions is a graph of codimension one in the phase
	space.

\vspace{6ex}

\noindent{\bf Mathematics Subject Classification 2010 (MSC2010):}
35K58, 35K90, 35B40.


\vspace{6ex}

\end{abstract}

 
\section{Introduction}

Let $\Omega\subseteq\re^{n}$ be a bounded connected domain.  In this
paper we consider the semilinear heat equation\begin{equation}
	u_{t}-\Delta u+|u|^{p}u=0
	\label{eqn:neumann}
\end{equation}
with homogeneous Neumann boundary conditions on $\partial \Omega$.  

The asymptotic behavior of solutions, and in particular their decay
rate and asymptotic profile as $t\to +\infty$, has been investigated
in the last decade by the third author and collaborators.  The
starting observation is that the Neumann Laplacian, namely the linear
operator associated to (\ref{eqn:neumann}), has a nontrivial kernel
consisting of all  constant functions.  This leads to the
coexistence of solutions with different decay rates.  In particular,
all nonzero solutions to (\ref{eqn:neumann}) are either \emph{fast
solutions} that decay exponentially, or \emph{slow solutions} with a
decay rate proportional to $t^{-1/p}$.  This is the so-called
null-slow-fast alternative, and was observed for the first time
in~\cite{ben-arbi} (see also~\cite{bah-sing}).  Similar results have
been obtained in~\cite{bah-d,hjk} for solutions to semilinear heat
equations such as $$ u_{t}-\Delta u+|u|^{p}u- \lambda_1 u=0$$ 
with homogeneous Dirichlet boundary conditions.  In the concrete model
(\ref{eqn:neumann}), as well as in the Dirichlet case, it has been
shown earlier in \cite{ben-arbi, bah-d,hjk} that all positive
solutions are slow.  In this paper we limit ourselves to the model
(\ref{eqn:neumann}), and we investigate more completely the sets of
initial data giving rise to slow/fast solutions.  We provide two
results.
\begin{enumerate}
	\renewcommand{\labelenumi}{(\arabic{enumi})} 
	
	\item In Theorem~\ref{prop:n-char-slow} we characterize fast
	solutions as those solutions which assume both positive and
	negative values for every $t\geq 0$.  What we actually prove is
	the contrapositive, namely that slow solutions are either
	eventually positive or eventually negative.  This motivates us to
	introduce the notion of \emph{positive-slow} and
	\emph{negative-slow} solutions.

	\item In Theorem~\ref{thm:main-neumann} we describe the set of
	initial data giving rise to slow/fast solutions.  We show that in
	the phase space $L^{2}(\Omega)$ there are two nonempty open sets
	of initial data originating positive-slow and negative
	slow-solutions, respectively.  These two open sets are separated
	by the graph of a continuous function, and all initial data in
	this graph give rise to a fast (or null) solution.
	
\end{enumerate}

Our characterization of slow/fast solutions in terms of their sign
follows from comparison principles and  the
construction of suitable subsolutions and supersolutions.

The characterization of Theorem~\ref{prop:n-char-slow} is the fundamental tool 
in the proof of Theorem~\ref{thm:main-neumann}. Indeed it implies 
that ``being the initial datum of a positive/negative slow solution'' 
is an open condition. At this point a simple connectedness argument 
implies that something different should exist in between, and by the 
null-slow-fast alternative the only remaining option is a fast or null 
solution.

\medskip

This paper is organized as follows.  In order to make the presentation
as self-contained as possible, in section~\ref{sec:well-known} we
collect all we need concerning existence, regularity and decay for
solutions to (\ref{eqn:neumann}).  In section~\ref{sec:statements} we
state our main results.  In section~\ref{sec:proof} we provide the proofs.
The final section~\ref{sec:rem} is devoted to some comments on
possible extensions of the results.


\setcounter{equation}{0}
\section{Basic tools and previous results}\label{sec:well-known}

Equation (\ref{eqn:neumann}) has been deeply investigated in
mathematical literature.  For the convenience of the reader, we
collect in this section the results that are needed in the sequel.

To begin with, we observe that (\ref{eqn:neumann}) can be interpreted
as the gradient flow in $L^{2}(\Omega)$ of the convex functional 
defined by
$$F_{p}(u):=\left\{
\begin{array}{ll}
	\displaystyle\frac{1}{2}\int_{\Omega}|\nabla u(x)|^{2}\,dx+
	\frac{1}{p+2}\int_{\Omega}|u(x)|^{p+2}\,dx\quad	 & 
	\mbox{if }u\in H^1(\Omega)\cap L^{p+2}(\Omega),  \\
	\noalign{\vspace{2ex}}
	+\infty & \mbox{otherwise}.
\end{array}
\right.
$$

Assuming that the boundary of $\Omega$ is $C^2$, the subdifferential
$B$ of $F_{p}$ is the operator defined by 
$Bu := -\Delta u + |u|^{p}u$ in the domain 
$$D(B) := \left\{u\in H^{2}(\Omega)\cap
L^{2p+2}(\Omega): \partial u/\partial n = 0
\mbox{ on }\partial \Omega \right\}.$$

As a consequence of \cite{OMM}, the operator $B$ generates a
contraction semigroup on $L^{2}(\Omega)$.  This provides existence,
uniqueness, and continuous dependence on initial conditions of a weak
solution $u \in C^{0}\left([0,+\infty),L^{2}(\Omega)\right)$ of
(\ref{eqn:neumann}) for any initial condition $u_0\in L^{2}(\Omega)$.
In the next statements we collect some well-known properties which
shall be used in the proofs of our main results.

\begin{thmbibl}[Regularity]\label{thmbibl:big} 
	
	Let $\Omega\subseteq\re^{n}$ be a bounded open set with boundary
	of class $C^{2}$, and let $p$ be a positive real number.  Let
	$u(t,x)$ be the unique solution to equation~(\ref{eqn:neumann})
	with homogeneous Neumann boundary conditions and initial datum
	$u_{0}\in L^{2}(\Omega)$, as defined previously.
	
	Then it turns out that
	\begin{equation}
		u\in W^{1, \infty}\left([\delta,+\infty),L^{2}(\Omega)\right)\cap
		L^{\infty}\left([\delta,+\infty),
		H^{2}(\Omega)\cap C(\overline{\Omega})\right)
		\quad\quad
		\forall\delta>0.
		\label{space:reg}
	\end{equation}
	
	Moreover, the solution satisfies the homogeneous Neumann boundary
	conditions in the classical sense for every $t>0$.
		
\end{thmbibl}

The second result concerns the comparison between two solutions with
different initial data.  From the general semigroup theory we know
that solutions depend continuously on initial data in $L^{2}(\Omega)$.
Here we need more, namely that the semigroup preserves the order, and
that one can estimate the norm in $L^{\infty}(\Omega)$ of the
difference between two solutions at positive times in terms of the
norm in $L^{2}(\Omega)$ of the difference between initial data.

\begin{thmbibl}[Comparison between two solutions]\label{thmbibl:comparison}
	Let $\Omega$ and $p$ be as in Theorem~\ref{thmbibl:big}. 
	Let $u(t,x)$ and $v(t,x)$ be two solutions to (\ref{eqn:neumann}) 
	with homogeneous Neumann boundary conditions and initial data 
	$u_{0}$ and $v_{0}$, respectively. 
	
	Then the following statements hold true.
	\begin{enumerate}
		\renewcommand{\labelenumi}{(\arabic{enumi})}
				
		\item \emph{(Order preservation)} If $u_{0}(x)\geq v_{0}(x)$
		for almost every $x\in\Omega$, then $u(t,x)\geq v(t,x)$ for
		every $t>0$ and every $x\in\overline{\Omega}$.
		
		\item \emph{(Continuous dependence $L^{2}(\Omega)\to 
		L^{\infty}(\Omega)$)} There exists a function
		$M_{1}:(0,+\infty)\to(0,+\infty)$ such that
		\begin{equation}
			|u(t,x)-v(t,x)|\leq
			M_{1}(t)\|u_{0}-v_{0}\|_{L^{2}(\Omega)}
			\quad\quad
			\forall t>0,\ \forall x\in\overline{\Omega}.
			\label{th:appl-l2-xd}
		\end{equation}
		
	\end{enumerate}
		
\end{thmbibl}

Although the estimate (\ref{th:appl-l2-xd}) is rather classical and
variants have been used in various contexts, for the sake of
completeness we provide a sketch of proof under the sole assumption
that a Sobolev-like imbedding $H^1(\Omega)\subseteq L^q(\Omega)$ is
satisfied for some $q>2$.

\begin{lemma}\label{lemma:moser} 
	Let us assume that there exists $q>2$ such that
	$H^1(\Omega)\subseteq L^q(\Omega) $ with continuous embedding, 
	namely there exists a constant $K_{0}$ such that
	\begin{equation}
		\|w\| _{ L^q(\Omega)} \le K_{0}\|w\|_{H^1(\Omega)}
		\quad\quad
		\forall w\in H^1(\Omega).
		\label{hp:moser-q}
	\end{equation}
	
	Let $T\in(0,1)$, let $c\in L^\infty((0,T)\times\Omega)$ be a 
	nonnegative function, and let
	$$z \in W^{1, \infty}\left((0,T),L^2(\Omega)\right)\cap 
	L^\infty\left((0,T),H^2(\Omega)\cap L^\infty(\Omega)\right)$$ 
	be a solution of
	\begin{equation}
		z_{t}-\Delta z+c(t,x)z=0
		\label{eqn:moser-z}
	\end{equation}
	with homogeneous Neumann boundary conditions.
	
	Then, setting $\beta:=q/(2q-4)$, it turns out that 
	\begin{equation}
		\| z(t,x) \|_{L^\infty(\Omega)} \leq
		\frac{4^{\beta^{2}}\cdot K_{0}^{2\beta}}{t^{\beta}}
		\|z(0,x)\|_{L^2(\Omega)}
		\quad\quad
		\forall t\in(0,T].
		\label{th:lemma-moser}
	\end{equation}

\end{lemma}

\paragraph{\textmd{\emph{Proof}}}

Let $r$ be any nonnegative real number. Let us multiply 
(\ref{eqn:moser-z}) by $|z|^{2r}z$ and let us integrate over 
$\Omega$. After integrating by parts the term $|z|^{2r}z\Delta z$, 
and recalling that $c(t,x)$ is nonnegative, we obtain that
\begin{equation}
	\frac{1}{2r+2}\frac{d}{dt}\int_{\Omega}|z|^{2r+2}\,dx+
	(2r+1)\int_{\Omega}|z|^{2r}|\nabla z|^{2}\,dx=
	-\int_{\Omega}c(t,x)|z|^{2r+2}\,dx\leq 0.
	\label{est:deriv-moser}
\end{equation}

This implies in particular that 
\begin{equation}
	\mbox{the function }t\to\|z(t,x)\|_{L^{\alpha}(\Omega)}\mbox{ is 
	nonincreasing for every $\alpha\geq 2$.}
	\label{est:monot-moser}
\end{equation}

Now we introduce the function $\psi_{r}(\sigma):=|\sigma|^{r}\sigma$, 
and we observe that 
$$\nabla[\psi_{r}(z)]=(r+1)|z|^{r}\nabla z.$$

As a consequence, (\ref{est:deriv-moser}) can be rewritten as
$$\frac{1}{2r+2}\frac{d}{dt}\left(
\|z\|_{L^{2r+2}(\Omega)}^{2r+2}\right)+
\frac{2r+1}{(r+1)^{2}}
\|\nabla[\psi_{r}(z)]\|_{L^{2}(\Omega)}^{2}\leq 0.$$

Given any $\tau\in(0,T]$, integrating in $[0,\tau]$ we deduce that
\begin{equation}
	\int_{0}^{\tau}
	\|\nabla[\psi_{r}(z(s,x))]\|_{L^{2}(\Omega)}^{2}\,ds\leq
	\frac{1}{2}\|z(0,x)\|_{L^{2r+2}(\Omega)}^{2r+2}.
	\label{est:psi-1}
\end{equation}

On the other hand, from (\ref{est:monot-moser}) with $\alpha:=2r+2$, 
we obtain also that
\begin{equation}
	\int_{0}^{\tau}
	\|\psi_{r}(z(s,x))\|_{L^{2}(\Omega)}^{2}\,ds\leq
	\tau\|z(0,x)\|_{L^{2r+2}(\Omega)}^{2r+2}.
	\label{est:psi-2}
\end{equation}

Adding (\ref{est:psi-1}) and (\ref{est:psi-2}) we conclude that
\begin{equation}
	\int_{0}^{\tau}
	\|\psi_{r}(z(s,x))\|_{H^{1}(\Omega)}^{2}\,ds\leq
	\left(\tau+\frac{1}{2}\right)
	\|z(0,x)\|_{L^{2r+2}(\Omega)}^{2r+2}
	\quad\quad
	\forall\tau\in(0,T].
	\label{est:moser-int}
\end{equation}

Now we exploit the continuous imbedding (\ref{hp:moser-q}). From 
(\ref{est:monot-moser}) with $\alpha:=(r+1)q$ we obtain that
\begin{eqnarray*}
	\tau\|\psi_{r}(z(\tau,x))\|_{L^{q}(\Omega)}^{2} & \leq & 
	\int_{0}^{\tau}
	\|\psi_{r}(z(s,x))\|_{L^{q}(\Omega)}^{2}\,ds    \\
	 & \leq & K_{0}^{2}\int_{0}^{\tau}
	 \|\psi_{r}(z(s,x))\|_{H^{1}(\Omega)}^{2}\,ds\\
	 & \leq & K_{0}^{2}\left(\tau+\frac{1}{2}\right)
	 \|z(0,x)\|_{L^{2r+2}(\Omega)}^{2r+2},
\end{eqnarray*}
and hence
$$\|\psi_{r}(z(\tau,x))\|_{L^{q}(\Omega)}^{2}\leq
K_{0}^{2}\left(1+\frac{1}{2\tau}\right)
\|z(0,x)\|_{L^{2r+2}(\Omega)}^{2r+2}
\quad\quad
\forall\tau\in(0,T].$$

Setting $\alpha:=2r+2$ and $\lambda:=q/2$, this can be written in the 
more suggestive form
$$\|z(\tau,x)\|_{L^{\lambda\alpha}(\Omega)}\leq
\left[K_{0}^{2}\left(1+\frac{1}{2\tau}\right)\right]^{1/\alpha}
\|z(0,x)\|_{L^{\alpha}(\Omega)}
\quad\quad
\forall\tau\in(0,T].$$

Due to the time-translation invariance, this implies also that
\begin{equation}
	\|z(\theta+\tau,x)\|_{L^{\lambda\alpha}(\Omega)}\leq
	\left[K_{0}^{2}\left(1+\frac{1}{2\tau}\right)\right]^{1/\alpha}
	\|z(\theta,x)\|_{L^{\alpha}(\Omega)}
	\label{est:moser-norm}
\end{equation}
whenever $0\leq\theta<\theta+\tau\leq T$.

This is the starting point of a classical iteration procedure. Given 
any $t\in(0,T]$, for every $n\in\n$ we set
$$t_{n}:=\left(1-\frac{1}{2^{n}}\right)t,
\hspace{4em}
\lambda_{n}:=2\lambda^{n},$$
and from (\ref{est:moser-norm}) with $\theta:=t_{n}$, 
$\tau:=t_{n+1}-t_{n}$, and $\alpha:=\lambda_{n}$ we deduce that
$$\|z(t_{n+1},x)\|_{L^{\lambda_{n+1}}(\Omega)}\leq
\left[K_{0}^{2}\left(1+\frac{2^{n}}{t}\right)\right]^{1/\lambda_{n}}
\|z(t_{n},x)\|_{L^{\lambda_{n}}(\Omega)}
\quad\quad
\forall n\in\n.$$

Since $t\leq 1$, this implies the simpler formula
$$\|z(t_{n+1},x)\|_{L^{\lambda_{n+1}}(\Omega)}\leq
\left[\frac{K_{0}^{2}\cdot 2^{n+1}}{t}\right]^{1/\lambda_{n}}
\|z(t_{n},x)\|_{L^{\lambda_{n}}(\Omega)}
\quad\quad
\forall n\in\n.$$

At this point an easy induction yields
\begin{equation}
	\|z(t_{n},x)\|_{L^{\lambda_{n}}(\Omega)}\leq
	2^{\gamma_{n}}\left[\frac{K_{0}^{2}}{t}\right]^{\beta_{n}}
	\|z(0,x)\|_{L^{2}(\Omega)}
	\quad\quad
	\forall n\in\n,
	\label{est:moser-n}
\end{equation}
where
$$\beta_{n}:=\sum_{k=0}^{n-1}\frac{1}{\lambda_{k}}\leq
\sum_{k=0}^{\infty}\frac{1}{2\lambda^{k}}=
\frac{1}{2}\frac{\lambda}{\lambda-1}=\frac{q}{2(q-2)}$$
and
$$\gamma_{n}:=\sum_{k=0}^{n-1}\frac{k+1}{\lambda_{k}}\leq
\sum_{k=0}^{\infty}\frac{k+1}{2\lambda^{k}}=
\frac{1}{2}\frac{\lambda^{2}}{(\lambda-1)^{2}}=
\frac{q^{2}}{2(q-2)^{2}}.$$

Letting $n\to+\infty$ in (\ref{est:moser-n}), we obtain 
(\ref{th:lemma-moser}).\qed

\bigskip
 
We are now ready to prove estimate (\ref{th:appl-l2-xd}). We consider
first the case where $u_{0}$ and $v_{0}$ are of class $C^2$ with compact
support in $\Omega$.  In this case $u$, $v$ and $z := u-v$ are bounded on
$(0,t) $ with values in $H^2(\Omega)\cap L^\infty(\Omega)$, and $z$ 
satisfies (\ref{eqn:moser-z}) with 
$$ c(t, x) : = \frac{|u|^{p}u- |v|^{p}v}{u-v}\ge 0 $$
with the convention that the quotient is $0$ whenever the denominator
vanishes.  Moreover, the continuous embedding $H^1(\Omega)\subseteq
L^q(\Omega)$ holds true with any $q\geq 2$ if $n\leq 2$, and with
$q=2^{*}:=2n/(n-2)$ if $n\geq 3$.  Therefore, we are in a position to
apply Lemma~\ref{lemma:moser}, from which we obtain
(\ref{th:appl-l2-xd}) with a function $M_{1}(t)$ independent of the
initial data.  At this point, the result for general initial data
follows from a density argument.\qed
 
\bigskip
 

The next statement describes all possible decay rates and
asymptotic profiles for solutions to (\ref{eqn:neumann}).  We refer
to~\cite[Theorem~1.3]{ben-arbi} and~\cite[Theorem~4.4]{ggh:casc-parab}
for further details and proofs.

\begin{thmbibl}[Classification of decay rates]\label{thmbibl:n-classification}
	
	Let $\Omega$ and $p$ be as in Theorem~\ref{thmbibl:big} and
	assume, in addition, that $\Omega$ is connected.  Let $u(t,x)$ be
	any solution to (\ref{eqn:neumann}) with homogeneous Neumann
	boundary conditions and initial datum in $L^{2}(\Omega)$. 
	
	Then
	one and only one of the following statements apply to $u(t,x)$.
	\begin{enumerate}
		\renewcommand{\labelenumi}{(\arabic{enumi})}
		
		\item \emph{(Null solution)} The solution is the null solution
		$u(t,x)\equiv 0$.
		
		\item \emph{(Slow solutions)} There exist $t_{0}>0$ and
		$M_{2}\geq 0$ such that
		$$\left||u(t,x)|-\frac{1}{(pt)^{1/p}}\right|\leq
		\frac{M_{2}}{t^{1+1/p}}
		\quad\quad
		\forall t\geq t_{0},\ \forall x\in\overline{\Omega}.$$
		
		\item \emph{(Spectral fast solutions)} There exist an
		eigenvalue $\lambda>0$ of the Neumann Laplacian, and a
		corresponding eigenfunction $\varphi_{\lambda}(x)$, such that
		$$\lim_{t\to+\infty}
		\left\|u(t,x)-\varphi_{\lambda}(x)e^{-\lambda t}\right\|_{L^{2}(\Omega)}
		e^{\gamma t}=0$$
		for some $\gamma>\lambda$.
	\end{enumerate}
	
\end{thmbibl}


\setcounter{equation}{0}
\section{Statements}\label{sec:statements}

In the first result of this paper we characterize slow solutions in
terms of sign.

\begin{thm}[Characterization of slow solutions]\label{prop:n-char-slow}
	Let $n$ be a positive integer, let $\Omega\subseteq\re^{n}$ be a
	bounded connected open set with $C^2$  boundary, and
	let $p$ be a positive real number.  Let $u(t,x)$ be a solution to
	equation (\ref{eqn:neumann}) with homogeneous Neumann boundary
	conditions.
	
	Then the following three statements are equivalent.
	\begin{enumerate}
		\renewcommand{\labelenumi}{(\roman{enumi})} 
		
		\item  There exist $t_{0}> 0$ and $M_{3}\geq 0$ such that 
		\begin{equation}
			\left||u(t,x)|-\frac{1}{(pt)^{1/p}}\right|\leq
			\frac{M_{3}}{t^{1+1/p}}
			\quad\quad
			\forall t\geq t_{0},\, \,\forall x\in\overline{\Omega}.			
			\label{th:n-slow}
		\end{equation}
	
		\item  There exist $t_{0}> 0$ and 
		$m:[t_{0},+\infty)\to(0,+\infty)$ such that
		\begin{equation}
			|u(t,x)|\geq m(t)
			\quad\quad
			\forall t\geq t_{0},\, \, \forall x\in\overline{\Omega}.
			\label{th:n-slow-2}
		\end{equation}
	
		\item There exists $t_{0}\geq 0$ such that either
		$u(t_{0},x)\geq 0$ for almost every $x\in\Omega$ or
		$u(t_{0},x)\leq 0$ for almost every $x\in\Omega$, but
		$u(t_{0},x)$ is not identically~0 in $\Omega$ (in the almost
		everywhere sense).
	\end{enumerate}
\end{thm}

Theorem~\ref{prop:n-char-slow} above implies that there are only two
types of slow solutions:
\begin{itemize}
	\item \emph{positive-slow solutions}, which are
	eventually positive and decay as $(pt)^{-1/p}$,

	\item \emph{negative-slow solutions}, which are eventually
	negative and decay as $-(pt)^{-1/p}$.
\end{itemize}

Moreover,  statement~(iii) implies that fast solutions are
necessarily sign changing functions for every $t\geq 0$.

\bigskip

The main result of this paper concerns the structure of slow/fast
solutions.  We show that, in the phase space $L^{2}(\Omega)$,
positive-slow and negative-slow solutions are separated by a manifold
of codimension one consisting of fast solutions.  As a consequence,
the set of initial data generating slow solutions is open and dense in
$L^{2}(\Omega)$.  The separating manifold is the graph of a continuous
function $\Phi$ defined in subspace $\np$ orthogonal to constant
functions (which are the kernel of the Neumann Laplacian).  The
function $\Phi$ turns out to be Lipschitz continuous when restricted
to $L^{\infty}(\Omega)$.

\begin{thm}[Structure of slow/fast solutions]\label{thm:main-neumann}
	
	Let $\Omega$ and $p$ be as in Theorem~\ref{prop:n-char-slow}.  Let
	us consider the space
	$$\np:=\left\{w\in L^{2}(\Omega):
	\int_{\Omega}w(x)\,dx=0\right\}.$$
	
	Then there exists a continuous function $\Phi:\np\to\re$ with the
	following property.  For every $w_{0}\in \np$, the solution
	$u(t,x)$ to equation (\ref{eqn:neumann}) with homogeneous Neumann
	boundary conditions and initial datum $u(0,x)=w_{0}(x)+k$ is
	\begin{itemize}
		\item  positive-slow if $k>\Phi(w_{0})$,
	
		\item  fast if $k=\Phi(w_{0})$ (or null if $w_{0}=0$, in 
		which case also $\Phi(w_{0})=0$),
	
		\item  negative-slow if $k<\Phi(w_{0})$.
	\end{itemize}
	
	Moreover, the function $\Phi$ is 1-Lipschitz continuous if 
	restricted to $L^{\infty}(\Omega)$, namely
	$$|\Phi(w_{1})-\Phi(w_{2})|\leq\|w_{1}-w_{2}\|_{L^{\infty}(\Omega)}
	\quad\quad
	\forall (w_{1},w_{2})\in
	\left[\np\cap L^{\infty}(\Omega)\right]^{2}.$$
\end{thm}

\begin{rmk}
	\begin{em}
		The projections of a function $u_{0}\in L^{2}(\Omega)$ on the 
		kernel of the Neumann Laplacian and on the orthogonal space 
		$\np$ are given, respectively, by
		$$\int_{\Omega}u_{0}(x)\,dx,
		\hspace{4em}
		u_{0}(x)-\int_{\Omega}u_{0}(x)\,dx.$$
		
		Therefore, the result of Theorem~\ref{thm:main-neumann} above
		can be rephrased by saying that the solution to
		(\ref{eqn:neumann}) with some initial datum $u_{0}\neq 0$ is
		positive-slow, fast, or negative-slow according to the sign of
		$$\int_{\Omega}u_{0}(x)\,dx-\Phi\left(
		u_{0}(x)-\int_{\Omega}u_{0}(x)\,dx\right).$$
		
	\end{em}
\end{rmk}



\setcounter{equation}{0}
\section{Proofs}\label{sec:proof}

To begin with, we recall that for every $p>0$ there exists a constant
$K_{p}$ such that
\begin{equation}
	\left||a+b|^{p}(a+b)-|a|^{p}a\strut\right|\leq
	K_{p}\left(|a|^{p}+|b|^{p}\right)|b|
	\quad\quad
	\forall(a,b)\in\re^{2}.
	\label{ineq:p-classic}
\end{equation}

This inequality follows from the mean value theorem applied to
the function $|x|^{p}x$.

\begin{lemma}\label{lemma:n-slow} 
	
	Let $n$, $\Omega$, and $p$ be as in
	Theorem~\ref{prop:n-char-slow}.  Let $v_{0}$ and $w_{0}$ be two
	functions in $L^{2}(\Omega)$ such that $v_{0}(x)\geq w_{0}(x)$ for
	almost every $x\in\Omega$, and $v_{0}(x)>w_{0}(x)$ on a set of
	positive measure.  Let $v$ and $w$ be the solutions to equation
	(\ref{eqn:neumann}) with homogeneous Neumann boundary conditions
	and initial data $v_{0}$ and $w_{0}$, respectively.
	
	Let us assume that $w$ is a fast or null solution in the sense of 
	Theorem~\ref{thmbibl:n-classification}.  
	
	Then $v$ is a slow solution in the sense of 
	Theorem~\ref{thmbibl:n-classification}.
	
\end{lemma}

\paragraph{\textmd{\emph{Proof}}}

Let $z(t,x):=v(t,x)-w(t,x)$ denote the difference, which is a 
nonnegative function because of statement~(1) of 
Theorem~\ref{thmbibl:comparison}, and satisfies
$$z_{t}=\Delta z-\left(|w+z|^{p}(w+z)-|w|^{p}w\strut\right).$$

Applying inequality (\ref{ineq:p-classic}) with $a:=w(t,x)$ and
$b:=z(t,x)$, we can estimate the nonlinear term in the right-hand
side, and obtain that 
$$z_{t}\geq \Delta z-K_{p}\left(|z|^{p}+|w|^{p}\right)z.$$

Let us consider now the function $I:[0,+\infty)\to[0,+\infty)$ 
defined by
$$I(t):=\int_{\Omega}z(t,x)\,dx
\quad\quad
\forall t\geq 0.$$

Our assumption on $v_{0}$ and $w_{0}$ implies that $I(0)>0$.  Since
the function $I(t)$ is continuous in $[0,+\infty)$ (due to the
continuity of $v$ and $w$ with values in $L^{2}(\Omega)$), there 
exists $\delta>0$ such that $I(\delta)>0$.

Now we argue by contradiction.  Let us assume that $v$, as well as
$w$, is not a slow solution.  By
Theorem~\ref{thmbibl:n-classification}, this implies that $z$ decays
exponentially to 0 in $L^{2}(\Omega)$, and hence also in
$L^{\infty}(\Omega)$ because of statement~(2) of
Theorem~\ref{thmbibl:comparison}, and therefore there exist constants
$\nu>0$ and $C_{\delta}>0$ such that $$z_{t}\geq \Delta
z-C_{\delta}e^{-\nu t}z \quad\quad
\forall (t,x)\in[\delta,+\infty)\times\Omega.$$

Integrating over $\Omega$ we find that
$$I'(t)\geq -C_{\delta}e^{-\nu t}I(t)
\quad\quad
\forall t\geq \delta,$$
and hence
$$I(t)\geq I(\delta)
\exp\left(-C_{\delta}\int_{\delta}^{+\infty}e^{-\nu s}\,ds\right)
\geq I(\delta)\exp(-C_{\delta}/\nu)
\quad\quad
\forall t\geq\delta.$$

This contradicts the fact that $z$ tends to 0 in $L^{2}(\Omega)$.\qed

\subsection*{Proof of Theorem~\ref{prop:n-char-slow}}

Implications $(i)\Rightarrow(ii)\Rightarrow(iii)$ are almost trivial.

As for $(iii)\Rightarrow(i)$, assuming
for instance the positive sign, it is enough to apply
Lemma~\ref{lemma:n-slow} with $ w(t,x):=0$ and
$v(t,x):=u(t+t_{0},x)$.\qed

\subsection*{Proof of Theorem~\ref{thm:main-neumann}}

In this proof we deal with many different initial conditions. For this 
reason we adopt the semigroup notation, namely we write 
$S_{t}(v_{0})$ or $[S_{t}(v_{0})](x)$ in order to denote the solution 
at time $t$ which has $v_{0}$ as initial condition.

\subsubsection*{Existence and uniqueness}

For every $w_{0}\in L^{2}(\Omega)$, let $\mathcal{K}^{+}(w_{0})$
denote the set of real numbers $k$ for which the solution with initial
datum $w_{0}(x)+k$ is positive-slow, and let $\mathcal{K}^{-}(w_{0})$
denote the set of real numbers $k$ for which the solution is
negative-slow.

The main point is proving that $\mathcal{K}^{+}(w_{0})$ and
$\mathcal{K}^{-}(w_{0})$ are, respectively, an open right half-line
and an open left half-line, and this two half-lines are separated by a
unique element.  When $w_{0}\in N^{\perp}$, this separator is the
value $\Phi(w_{0})$ that we are looking for.  We prove these claims
through several steps.

\subparagraph{\textmd{\textit{Step 1}}}

We prove that $\mathcal{K}^{+}(w_{0})$ and $\mathcal{K}^{-}(w_{0})$
are nonempty for every $w_{0}\in L^{2}(\Omega)$.

To this end, we concentrate on $\mathcal{K}^{+}(w_{0})$, since the
argument for $\mathcal{K}^{-}(w_{0})$ is symmetric.  Let us choose two
positive constants $m_{0}$ and $t_{0}$, and let us set
$$m_{1}:=\frac{m_{0}}{2\left(1+pm_{0}^{p}t_{0}\right)^{1/p}},
\hspace{4em}
\delta:=\frac{m_{1}}{M_{1}(t_{0})},$$
where $M_{1}(t_{0})$ is the constant which appears in inequality
(\ref{th:appl-l2-xd}).

Let us choose $v_{0}\in C^{0}(\overline{\Omega})$ such that
$\|w_{0}-v_{0}\|_{L^{2}(\Omega)}\leq\delta$ (this is possible because
$C^{0}(\overline{\Omega})$ is dense in $L^{2}(\Omega)$).  Due to
boundedness of $v_{0}$, there exists $k_{0}\in\re$ such that
$$v_{0}(x)+k_{0}\geq m_{0}
\quad\quad
\forall x\in\overline{\Omega}.$$

Now we claim that
$$\left[S_{t}(v_{0}+k_{0})\right](x)\geq
\frac{m_{0}}{\left(1+pm_{0}^{p}t\right)^{1/p}}
\quad\quad
\forall t\geq 0,\ \forall x\in\overline{\Omega}.$$

This inequality follows from the usual comparison principle because it
is true when $t=0$, and in addition both the left-hand and the
right-hand side are solutions to (\ref{eqn:neumann}) with homogeneous
Neumann boundary conditions.  Setting $t=t_{0}$, and recalling our
definition of $m_{1}$, we obtain that
\begin{equation}
	\left[S_{t_{0}}(v_{0}+k_{0})\right](x)\geq 2m_{1}
	\quad\quad
	\forall x\in\overline{\Omega}.
	\label{est:d1}
\end{equation}

On the other hand, statement~(2) of Theorem~\ref{thmbibl:comparison}
applied to initial data $w_{0}+k_{0}$ and $v_{0}+k_{0}$
implies that
\begin{equation}
	|\left[S_{t_{0}}(w_{0}+k_{0})\right](x)-
	\left[S_{t_{0}}(v_{0}+k_{0})\right](x)|\leq
	M_{1}(t_{0})\|w_{0}-v_{0}\|_{L^{2}(\Omega)}\leq
	m_{1}
	\label{est:d2}
\end{equation}
for every $x\in\overline{\Omega}$.  From (\ref{est:d1}) and
(\ref{est:d2}) it follows that
$$\left[S_{t_{0}}(w_{0}+k_{0})\right](x)\geq m_{1} \quad\quad
\forall x\in\overline{\Omega}.$$

Thanks to Theorem~\ref{prop:n-char-slow}, this is enough to
conclude that the solution with initial condition
$w_{0}(x)+k_{0}$ is positive-slow, and hence
$k_{0}\in\mathcal{K}^{+}(w_{0})$.

\subparagraph{\textmd{\textit{Step 2}}}

We prove that $\mathcal{K}^{+}(w_{0})$ is an open right half-line, and
analogously $\mathcal{K}^{-}(w_{0})$ is an open left half-line.

Let us consider $\mathcal{K}^{+}(w_{0})$ (the argument for
$\mathcal{K}^{-}(w_{0})$ is symmetric).  It is a right half-line
because, if $w_{0}(x)+k_{0}$ gives rise to a positive-slow solution
$u(t,x)$, then every solution with initial datum $w_{0}(x)+k$ with
$k>k_{0}$ is greater than $u(t,x)$, and hence it is positive-slow as
well.

It remains to show that $\mathcal{K}^{+}(w_{0})$ is an open set.  Let
us assume that $k_{0}\in\mathcal{K}^{+}(w_{0})$, so that the solution
$u(t,x)$ with initial datum $w_{0}(x)+k_{0}$ is positive-slow.
Due to Theorem~\ref{prop:n-char-slow}, it turns out that
$$\left[S_{t_{0}}(w_{0}+k_{0})\right](x)\geq
m_{0}
\quad\quad
\forall x\in\overline{\Omega}$$
for suitable constants $t_{0}> 0$ and $m_{0}>0$.  Applying
statement~(2) of Theorem~\ref{thmbibl:comparison} as in the previous step, we
obtain that
$$\left[S_{t_{0}}(w_{0}+k)\right](x)\geq
\frac{m_{0}}{2}
\quad\quad
\forall x\in\overline{\Omega}$$
provided that $k$ is close enough to $k_{0}$. Applying
Theorem~\ref{prop:n-char-slow} once again, we can conclude that
all these neighboring solutions are positive-slow as well.

\subparagraph{\textmd{\textit{Step 3}}}

The structure of $\mathcal{K}^{+}(w_{0})$ and $\mathcal{K}^{-}(w_{0})$
implies that
\begin{equation}
	\sup\mathcal{K}^{-}(w_{0})\leq\inf\mathcal{K}^{+}(w_{0}),
	\label{sup-inf}
\end{equation}
and any $k$ in between (endpoints included) lies neither in
$\mathcal{K}^{-}(w_{0})$ nor in $\mathcal{K}^{+}(w_{0})$.  Due to the
null-slow-fast alternative of Theorem~\ref{thmbibl:n-classification},
the corresponding solutions are necessarily fast or null.  Finally,
as a consequence of Lemma ~\ref{lemma:n-slow}, we do have equality in
(\ref{sup-inf}), and hence for every $w_{0}\in L^{2}(\Omega)$ there 
exists a unique $k$ such that $w_{0}(x)+k$ generates a fast (or null) 
solution.

\subsubsection*{Continuity}

We show that the map $w_{0}\to\Phi(w_{0})$ is continuous with respect
to the norm of $L^{2}(\Omega)$, namely for every $w_{0}\in \np$
and every $\ep>0$ there exists $\delta>0$ such that
\begin{equation}
	\Phi(w_{0})-\ep\leq\Phi(v_{0})\leq\Phi(w_{0})+\ep
	\label{est:Phi-cont}
\end{equation}
for every $v_{0}\in \np$ with 
$\|v_{0}-w_{0}\|_{L^{2}(\Omega)}\leq\delta$.

Let us consider the solution with initial condition 
$w_{0}(x)+(\Phi(w_{0})+\ep)$. It is positive-slow, and hence 
from Theorem~\ref{prop:n-char-slow} we know that
$$\left[S_{t_{0}}(w_{0}+(\Phi(w_{0})+\ep))\right](x)\geq
m_{0}
\quad\quad
\forall x\in\overline{\Omega}$$
for suitable constants $t_{0}> 0$ and $m_{0}>0$.  Applying
statement~(2) of Theorem~\ref{thmbibl:comparison} as in the existence
part, we deduce that
$$\left[S_{t_{0}}(v_{0}+(\Phi(w_{0})+\ep))\right](x)\geq
\frac{m_{0}}{2}
\quad\quad
\forall x\in\overline{\Omega}$$
provided that $\|v_{0}-w_{0}\|_{L^{2}(\Omega)}$ is small enough.
Applying Theorem~\ref{prop:n-char-slow} once again, we deduce that
the solution with initial condition
$v_{0}(x)+(\Phi(w_{0})+\ep)$ is positive-slow as well. It follows that
$\Phi(w_{0})+\ep\in\mathcal{K}^{+}(v_{0})$, and therefore
$\Phi(v_{0})\leq\Phi(w_{0})+\ep$.

This proves that the inequality on the right in (\ref{est:Phi-cont})
holds true for every $v_{0}\in \np$ which is close enough to
$w_{0}$ with respect to the norm of $L^{2}(\Omega)$.  A symmetric
argument applies to the inequality on the left.

\subsubsection*{Lipschitz continuity with respect to the uniform norm}

We show that the map $w_{1}\to\Phi(w_{1})$ restricted to
$L^{\infty}(\Omega)$ is Lipschitz continuous with Lipschitz constant
equal to~1.  To this end, we take any $w_{1}$ and $w_{2}$ in $\np\cap
L^{\infty}(\Omega)$, and from the definition of $L^{\infty}(\Omega)$
we deduce that 
$$w_{1}(x)+(\Phi(w_{1})+\ep)\leq w_{2}(x)+
\left(\|w_{1}-w_{2}\|_{L^{\infty}(\Omega)}+\Phi(w_{1})+\ep\right)$$
for every $\ep>0$ and almost every $x\in\Omega$.  Since the solution
with initial datum equal to the left-hand side is positive-slow, the
solution with initial datum equal to the right-hand side is
positive-slow as well, which proves that
$$\|w_{1}-w_{2}\|_{L^{\infty}(\Omega)}+\Phi(w_{1})+
\ep\in\mathcal{K^{+}}(w_{2}).$$

In an analogous way we obtain that 
$$-\|w_{1}-w_{2}\|_{L^{\infty}(\Omega)}+\Phi(w_{1})-
\ep\in\mathcal{K^{-}}(w_{2}).$$

Since $\Phi(w_{2})$ separates $\mathcal{K^{-}}(w_{2})$ and
$\mathcal{K^{+}}(w_{2})$, letting $\ep\to 0^{+}$ we conclude
that 
$$\Phi(w_{1})-\|w_{1}-w_{2}\|_{L^{\infty}(\Omega)}
\leq\Phi(w_{2})\leq
\Phi(w_{1})+\|w_{1}-w_{2}\|_{L^{\infty}(\Omega)},$$
which completes the proof.\qed

\section{Additional results and possible extensions}\label{sec:rem} 

In this section we describe some additional properties and research
directions.

\paragraph{Strong positivity}

Under the $C^2$ regularity assumption on $\partial\Omega$ which
allowed us to set properly the problem, more can be said about the
behavior of solutions.  Actually, an application of the strong minimum
principle gives that for any non-negative initial value $u_0\in
L^2(\Omega)$ which is positive on a set of positive measure, the
solution of (\ref{eqn:neumann}) is uniformly (with respect to $x\in
\Omega $) positive for all positive times.  This of course means that
right after the first time at which a slow solution has a constant
sign, it becomes strictly above a positive (time depending) constant
or strictly below a time depending negative constant.

\paragraph{Relaxed regularity } 

In principle, in order for (\ref{eqn:neumann}) to be properly set in a
reasonable sense (for instance variational or distributional), a $C^1$
or even Lipschitz regularity assumption on $\partial\Omega$ seems to
be enough.  In such a case, the Sobolev embedding theorem would be
applicable, the difficulty might be that the solution does not need to
be continuous up to the boundary for $t>0$.  The relevant regularity
class for solutions would then be 
$$W^{1,\infty}\left([\delta,+\infty),L^{2}(\Omega)\right)\cap
L^{\infty}\left([\delta,+\infty), H^{1}(\Omega)\cap L^\infty
({\Omega})\cap C({\Omega})\right),$$
since interior regularity is always true.  This is enough to state
properly, mutatis mutandis, our various results.

\paragraph{More general nonlinearities}
		
For the sake of simplicity, we presented our results for the equation
with the model nonlinearity $|u|^{p}u$.  Nevertheless, the theory can
be extended with little effort to more general nonlinear terms
$f(u(t))$.  The essential assumption here is that $f:\re\to\re$ is an
increasing function such that $f(x)\sim|x|^{p}x$ and
$f'(x)\sim(p+1)|x|^{p}$ as $x\to 0$.  One can even assume monotonicity
just in a neighborhood of the origin, but of course in that case one
obtains a description of solutions only for initial data whose norm in
$L^{\infty}(\Omega)$ is small enough.
	
\paragraph{The Dirichlet case} 

Many results for concrete models have been unified
in~\cite{ggh:casc-parab} by developing an abstract theory for
evolution inequalities of the form
\begin{equation}
	|u'(t)+Au(t)|_{H}\leq
	K_{0}\left(|u(t)|_{H}^{1+p}+|A^{1/2}u(t)|_{H}^{1+q}\right)
	\quad\quad \forall t\geq 0,
	\label{pbm:diff-ineq}
\end{equation}
where $A$ is a self-adjoint nonnegative operator with discrete
spectrum in a Hilbert space $H$, and $K_{0}$, $p$, $q$ are positive
real numbers.  A full description of possible decay rates was
provided, showing that all nonzero solutions to (\ref{pbm:diff-ineq})
that decay to~0 are either exponentially fast as solutions to the
linearized equation $u'(t)+Au(t)=0$, or slow as solutions to the
ordinary differential inequality $|u'(t)|\leq K_{0}|u(t)|^{1+p}$.  By
relying on this kind of general techniques, it will be possible to
extend this theory to more general parabolic partial differential
equations whose linear part has a nontrivial kernel, for instance the
problem ``at resonance''
$$ u_{t}-\Delta u+|u|^{p}u- \lambda_1 u=0$$
with homogeneous Dirichlet boundary conditions.  In the present paper
we decided to limit ourselves to the model example (\ref{eqn:neumann})
in which some of the arguments appear simpler.  The other cases will
be studied elsewhere.

\paragraph{Second order equations} It might be interesting to look
for a concrete realization of the slow-fast alternative for second
order evolution equations with dissipative terms, in the same way as
the results of~\cite{ggh:casc-parab} were extended
in~\cite{ggh:jems,ggh:casc-hyperb}.  However this would require
completely new ideas, since both the regularizing effect and the
comparison principles are specific to parabolic problems and even in
the simple case of the ordinary differential equation 
$$ u''+u'+u^3 = 0 $$ 
the set of fast solutions has already a rather complicated shape.
Moreover, in that case, all non-trivial solutions (including
exponentially decaying ones) are asymptotically signed, so that even
for the hyperbolic problem $$ u_{tt}+u_{t}-\Delta u +u^3 = 0 $$ with
Neumann homogeneous boundary conditions, the slow character is not
equivalent to the existence of a constant sign for $t$ large.


\label{NumeroPagine}

\end{document}